\documentclass[12pt]{amsart}
\usepackage{xcolor}
\usepackage{graphicx}
\theoremstyle{definition}
\newtheorem{definition}{Definition}
\theoremstyle{plain}
\newtheorem{theorem}{Theorem}
\newtheorem{corollary}{Corollary}
\newtheorem{lem}{Lemma}
\newtheorem{prop}{Proposition}
\newtheorem{ejem}{Example}
\theoremstyle{remark}

\newtheorem{remark}{Remark}
\DeclareMathOperator{\CC}{\mathbb C}

\DeclareMathOperator{\NN}{\mathbb N}

\def\limsup{\mathop{\overline{\rm lim}}}

\title{BOUNDEDNESS  OF ZEROS OF SOBOLEV ORTHOGONAL POLYNOMIALS VIA GENERALIZED EIGENVALUES
 }
\author[1,2]{ Escribano, C.$^{1,2}$}

\author[1]{ Gonzalo, R.$^{1}$}

\address{[1] Departamento de Matem\'atica Aplicada a las Tecnolog\'{i}as de la Informaci\'on y las Comunicaciones, Escuela T\'ecnica Superior de Ingenieros Inform\'aticos y Center for Computational Simulation, Universidad Polit\'ecnica de Madrid, Campus de Montegancedo\\
Boadilla del Monte, Madrid Spain. \\
[2] Center for Computational Simulation, Universidad Polit\'ecnica de Madrid,Spain.\\
}

\email{cescribano@fi.upm.es}
\email{rngonzalo@fi.upm.es}

\thanks{The authors would like to thank the project Madrid Quantum-CM, funded by the European Union,  NextGenerationEU (PRTR-C17.I1) and by the Comunidad de Madrid, Programa de Acciones Complementarias
}

\begin{document}

\begin{abstract} The main aim of this work is to apply the study of the asymptotic behaviour of generalized eigenvalues between infinite Hermitian definite positive matrices in an important question regarding the  location of zeros of Sobolev orthogonal polynomials.  In order to do it we introduce  matrix  Sobolev inner products  associated with a set  of infinite Hermitian positive definite matrices that generalize a type of Sobolev inner products. These matrix components are not  necessarily moment matrices associated to measures. This  general framework   allows us  to  study boundedness of  multiplication operator on the space of polynomials, that is a sufficient condition for boundedness of zeros of orthogonal polynomials,   via this matrix approach. As a consequence of our results we provide a criteria  to assure boundedness of the zero set of Sobolev orthogonal polynomials in different situations in terms of the generalized eigenvalues introduced in \cite{EGT6}.

\end{abstract}

\maketitle
\begin{quotation} {\sc {\footnotesize Keywords}}. {\small Hermitian matrices, matrix Sobolev inner product, zeros of orthogonal  polynomials, generalized eigenvalue, measures.}
\end{quotation}

\section{Introduction}

\noindent  Let $\mathbf{M}=(c_{i,j})_{i,j=0}^{\infty}$ be an infinite Hermitian matrix, i.e., $c_{i,j}= \overline{c_{j,i}}$ for all
$i,j$ non-negative integers. Following \cite{EGT1} we say that $\mathbf{M}$ is positive definite (resp. positive semidefinite)  if
all the truncated sections are positive definite (resp. positive semidefinite). An
Hermitian positive definite  matrix (in short, an HPD matrix) defines an inner product $\langle \cdot, \cdot \rangle_{\mathbf{M}}$ in the space $\mathbb{P}[z]$ of all polynomials with complex coefficients in the following way:
for   $p(z)=\sum_{k=0}^{n}v_kz^k, 
q(z)=\sum_{k=0}^{m}w_kz^k \in \mathbb{P}[z]$ 

\begin{equation}\label{aste}
    \langle p(z)
,q(z) \rangle_{\mathbf{M}}=v\mathbf{M}w^{*},
\end{equation}
\\
\noindent being $v=(v_0,\dots,v_n,0,0, \dots), w=(w_0,\dots,w_m,0,0, \dots) \in c_{00}$ where $c_{00}$ is the space
of all complex sequences with only finitely many non-zero entries. In the case that $\mathbf{M}$ is positive definite this inner product induces an Hilbertian norm 
denoted by  $\Vert p(z) \Vert_{\mathbf{M}}^2 = \langle p(z),p(z)\rangle_{\mathbf{M}} $ for every $p(z)\in \mathbb{P}[z]$ and the completion of polynomials with such a norm, denoted by  $P^{2}(\mathbf{M})$, is an infinite dimensional  Hilbert space.

\bigskip

\noindent An  interesting
class of infinite Hermitian semidefinite positive matrices (in short HSPD matrices) are those which are  moment matrices
 with respect to a certain positive 
 measure $\mu$, i.e.,  $\mathbf{M}=(c_{n,m})_{n,m=0}^{\infty}$ such that there exists a {\it representating} measure $\mu$ with support on $\CC$ and  finite moments for all $n,m\geq 0$ verifying 
$$
c_{n,m}= \int z^{n} \overline{z}^md\mu.
$$
\noindent Thorough the paper we  always consider Borel positive compactly supported measures in the complex plane and with finite moments. The associated moment matrix to a measure $\mu$ (denoted by $\mathbf{M}(\mu)$) is an HPD matrix if and only if the support of $\mu$ is an infinite set. 

\noindent Associated to an infinite HPD matrix, not necessarily a moment matrix,  and via Gram-Schmidt procedure  we may construct  a sequence of orthonormal polynomials, which is unique if the leading  coefficients are positive.  The essential key to relate  the theory of orthogonal polynomials with  matrix analysis of infinite matrices, via the  moment matrix as in  \cite{EGT1}, is the following identity: being $p(z)=v_0+v_1z+\dots+v_nz^n \in \mathbb{P}[z], $  $(v_0,v_1,\dots,v_n,0,0,0,\dots)\in c_{00}$,

\begin{equation}\label{aste2}
    v\mathbf{M}(\mu)v^{*}=\int \vert p(z) \vert^2 \; d\mu \qquad  \qquad  
\end{equation}

\noindent Other important class of  HPD
matrices are Sobolev moment matrices with respect to a set of Borel measures  $\{\mu_{j}\}_{j=0}^{k}$  compactly supported in the complex plane and such that  the support of $\mu_0$ contains infinitely many points \cite{Finkelshtein}. The inner product in $\mathbb{P}[z]$ is defined as follows,
$$
<p(z),q(z)>_{\mathbf{S}}=\sum_{j=0}^{k} \int p^{(j)}(z)\overline{q^{(j)}(z)}d\mu_k(z), \qquad p(z), q(z)\in \mathbb{P}[z]
$$
\noindent where as usual $p^{(j)}$ denotes the $j$-th derivative of a polynomial. The Sobolev moment matrix associated to this inner product is   $\mathbf{S}:=(c_{n,m})_{n,m=0}^{\infty}$ where for $n,m\geq 0$
$$
c_{n,m}=<z^{n},z^{m}>_{\mathbf{S}}.
$$
\noindent   As it is observed in \cite{Finkelshtein} for  Hankel matrices, and can be easily generalized to the more general context of HPD matrices, denoting by  $\mathbf{M}_j:=\mathbf{M}(\mu_j)$ the associated moment matrix with $\mu_j$ for $j=0,1\dots k$, the Sobolev moment matrix $\mathbf{M}_{\mathbf{S}}$ is given by
$$
\mathbf{M}_{\mathbf{S}}=\sum_{j=0}^{k} \mathbf{\Lambda}^{j}\mathbf{M}_j\left(\mathbf{\Lambda}^{j}\right)^{*}
$$

\noindent where
\begin{equation}
\mathbf{\Lambda}^{j}=(\lambda^{j}_{n,m})_{n,m=0}^{\infty}, \qquad \lambda^{j}_{n,m}=\dfrac{n!}{(n-j)!}\delta_{j,n-m}, j=0,\dots,k.
\end{equation} \label{Filk}

\noindent In an analogous way, following the matrix approach, we here define  {\it matrix Sobolev inner products} associated to a certain  set of  HPSD matrices $\{\mathbf {M}_j\}_{k=0}^{k}$,  with $\mathbf{M}_0$ being HPD in $\mathbb{P}[z]$ as

$$<p(z),q(z)>_{\mathbf{M}_{S}}=v\sum_{j=0}^{k}\mathbf{\Lambda}^{j}\mathbf{M}_j\left(\mathbf{\Lambda}^{j}\right)^{*} w^{*}$$   
where ${\mathbf{M}_{S}}=\sum_{j=0}^{k} \mathbf{\Lambda}^{j}\mathbf{M}_j\left(\mathbf{\Lambda}^{j}\right)^{*} $ and $\mathbf{{\Lambda}}^{j}$ as above. This is a more general framework than  Sobolev inner products on the space of polynomials.

\bigskip
\noindent In this paper we are concerned with boundedness  of  zeros of orthogonal polynomials associated to a matrix Sobolev inner product using the matrix approach. Following this matrix via  
we continue the study of the behaviour of the smallest and largest generalized eigenvalues initiated in  \cite{EGT6} in relation with this type of problems. 

\bigskip

\noindent The  location of the zeros of orthogonal  polynomials  is a problem of vital importance to study the asymptotics of orthogonal polynomials (see e.g. \cite{LP},\cite{LPP}). Although the boundedness of the set of zeros of orthogonal polynomials holds for every measure with  compact  support
in the complex plane, it is an open problem to bound the zeros of Sobolev orthogonal polynomials associated to compactly supported measures   in the complex plane. The boundedness of the zeros is a consequence of the boundedness of the multiplication operator in the corresponding space of polynomials with  the associated  Sobolev inner product; moreover, the zeros of Sobolev orthogonal polynomials are contained in the disk centered at zero and whose  radio is the norm of the multiplication operator (see \cite{LPP}). In \cite{APRR}, there are some answers to the question stated in \cite{LP} about some conditions for the boundedness of multiplication operator. 
Indeed,  in \cite{APRR} they give a necessary and sufficient condition so that the multiplication operator associated to a vectorial measure is bounded in terms the existence of an extended sequentially dominated vector measure whose norms in the space of the polynomials $\mathbb{P}[z]$ are equivalent. 

\bigskip
\noindent In relation with the study of the asymptotic behaviour of smallest eigenvalues of  truncated moment matrices  in the theory of orthogonal polynomials, it was first introduced in \cite{Berg} (see also \cite{Berg-Duran}) regarding  the problem of the determinacy of a measure in the real line. Later on, for compactly supported measures in the complex plane in \cite{EGT1} this asymptotic behaviour  is studied in the problem of completeness of a measure $\mu$ in the sense of density of polynomials in the corresponding Hilbert  spaces $L^{2}(\mu)$ (see also \cite{EGT5}). Recently, in \cite{EGT6} the authors introduce the notion of smallest and largest generalized eigenvalues between two moment matrices in order to compare two different measures with support in the complex plane; such technique  provides generalizations of  the results in \cite{EGT5}. The study of generalized eigenvalues between moment matrices in \cite{EGT6}  is also  used in order to compare the supports of two measures.

\bigskip

\bigskip

\noindent The structure of the paper is as follows.

\noindent In  section 1  we prove that in the general framework of  inner products induced by   HPD matrices,   boundedness of the multiplication operator is a sufficient condition to assure the boundedness of the set of zeros of orthogonal polynomials. The proof is analogous to  the case of Sobolev inner products (see \cite{LP}, also\cite{Duran-Saff}). We also provide a sufficient condition to  boundedness of the multiplication operator in terms of certain properties of the sequence of entries of the HPD matrix. 

\medskip 
\noindent In section 2, motivated by the  notions of  sequentially dominated  and essentially sequentially dominated vectorial measures ( \cite{LP}, \cite{JR3}) and using a matrix approach we introduce the  notion   of {\it sequentially dominated} HPSD matrices. In this sense, we introduce the class of  {\it matrix sequentially dominated }(in short, MSD) vectorial measures as those that their associated moment matrices are sequentially dominated.  In particular,  MSD  vectorial measures are sequentially dominated; nevertheless,  the converse is not true as we show in Example 1.  

\bigskip
\noindent In the general  framework of matrix inner Sobolev products our   main result is  Theorem 1, where  we prove that the multiplication operator associated to a matrix Sobolev inner product of a sequentially dominated HPD matrices, $\{\mathbf{M}_j\}_{j=0}^{k}$, is bounded,  assuming that for all $j=0,\dots,k$  the corresponding  multiplication  operator, $\mathcal{D}_j$, associated  to $\mathbf{M}_j$  is bounded. Consequently, for this inner product the set of zeros of associated orthogonal polynomials is bounded. Moreover, by introducing a type of  comparison between matrix Sobolev inner products we prove that MSD vectorial  measures introduced here are, in particular,essentially sequentially dominated vectorial measures. 

\bigskip

\noindent Finally, in  the last section we characterize {\it sequentially dominated}  HPD  using the notion of largest generalized eigenvalues between infinite matrices. This approach lets  to apply the behaviour of these generalized eigenvalues in order to give a criteria  for boundedness of multiplication operator, and consequently, for boundedness of the set of zeros in some examples of Sobolev inner products. In particular, for the case of Sobolev inner products associated to two measures being  one of them  the Lebesgue measure in the unit circle we provide sufficient conditions in terms of the  asymptotic behaviour of the standard 
 smallest and largest eigenvalues for boundedness of  zeros of Sobolev orthogonal polynomials.

\bigskip 

\noindent We begin with some definitions:

\begin{definition} Let $\mathbf{A},\mathbf{B}$  be two infinite Hermitian matrices. We say that
\begin{enumerate} 
\item $\mathbf{B} \leq \mathbf{A}$ if
$
v\mathbf{B}v^{*} \leq v \mathbf{A}v^{*},\qquad  \text{ for every } v\in c_{00}$

\smallskip

\item 
$\mathbf{B} < \mathbf{A}$ if
$
v\mathbf{B}v^{*} < v \mathbf{A}v^{*},\qquad  \text{ for every } 0\neq v\in c_{00}$
\end{enumerate}
\end{definition}

\bigskip
\noindent In particular, for an hermitian matrix $\mathbf{A}$ we denote $\mathbf{A}>0$ (resp. $\mathbf{A} \geq 0$) if $\mathbf{A}$ is positive definite (resp. positive semidefinite).

\begin{definition} Let $\mathbf{A},\mathbf{B}$ be HPD matrices. We say that $\mathbf{A}$ is {\it comparable} with $\mathbf{B}$, and we denote it by $\mathbf{A} \asymp \mathbf{B}$, if there exist positive constants $c,C>0$ such that 
$$
c\mathbf{A} \leq \mathbf{B} \leq C \mathbf{A}
$$
\end{definition}


\begin{definition}Given two finite matrices of size $n\times n$,  $\mathbf{A}$ 
and $\mathbf{B}$,  we say that $\lambda$ is an eigenvalue of $\mathbf{A}$ with respect to $\mathbf{B}$ if there exists a nonzero vector $v\in \CC^{n}$ such that $\mathbf{A}v= \lambda \mathbf{B}v$. These are called the {\it generalized eigenvalues} of $\mathbf{A}$ with respect $\mathbf{B}$ .
\end{definition}

\noindent It is well known that in the
 case of hermitian matrices the generalized eigenvalues are real numbers. Note that in the case that the second matrix is the identity matrix we obtain the standard  eigenvalues. In \cite{EGT6} the following indexes are defined as a generalization of smallest and largest generalized eigenvalues for infinite dimensional matrices: 

\begin{definition} Let $\mathbf{A},\mathbf{B}$ be Hermitian matrices of size $n$, we define $\lambda_n(\mathbf{A},\mathbf{B})$ as the smallest eigenvalue of $\mathbf{A}$ with respect to $\mathbf{B}$ and
$\beta_n(\mathbf{A},\mathbf{B})$ as the largest eigenvalue of $\mathbf{A}$ with respect to $\mathbf{B}$.
\end{definition}

\begin{definition} Let  $\mathbf{A},\mathbf{B}$ two infinite HPD matrices. We define:
\begin{eqnarray*}
0 \leq \lambda(\mathbf{A},\mathbf{B})&:=&\lim_{n\to \infty} \lambda_{n}(\mathbf{A}_n,\mathbf{B}_n)\\
\beta(\mathbf{A},\mathbf{B})& := &  \lim_{n\to
\infty}\beta_{n}(\mathbf{A}_n,\mathbf{B}_n) \leq \infty
\end{eqnarray*}
\noindent where $\mathbf{A}_n, \mathbf{B}_n$ are the corresponding truncated matrices of size  $(n+1)\times (n+1)$ of $\mathbf{A},\mathbf{B}$ respectively.

\end{definition}


\section{Boundedness of multiplication operator and zeros of orthogonal polynomials}

\noindent Let $\mathbf{M}$ be an HPD matrix and let $P^{2}(\mathbf{M})$ be the Hilbert space endowed with the norm induced by $\mathbf{M}$. The multiplication operator is the linear mapping given by 
$$
\begin{array}{ccc}
   \mathcal{D}:\mathbb{P}[z] & \rightarrow & \mathbb{P}[z] \\
   p(z) &\rightarrow & zp(z)
 \end{array}
$$ 

\noindent  The operator $\mathcal{D}$ is bounded in $\mathbb{P}[z]$ if 
$$
\Vert \mathcal{D} \Vert = \sup \{ \Vert zp(z)\Vert_{\mathbf{M}}: \;  p(z)\in \mathbb{P}[z], \Vert p(z) \Vert_{\mathbf{M}}=1\}<\infty,
$$
\noindent and in this case  $\mathcal{D}$ can be extended to an operator in $P^{2}(\mathbf{M})$. 

\bigskip
\noindent We first prove that  boundedness of the multiplication operator on $\mathbb{P}[z]$ is a sufficient condition for the boundedness of the set of zeros of orthogonal polynomials associated with respect any  inner product induced by any HPD matrix. This  is  a generalization of the result for standard  orthogonal polynomials associated with a measure, for  the corresponding one  for  Sobolev orthogonal polynomials  given in \cite{LPP} and for nonstandard polynomials (\cite{Duran-Saff}). The proof is the same and we  include it for its simplicity and the sake of completeness. 

\begin{prop} Let  $\mathbf{M}$ be an infinite HPD matrix and let  $\{\varphi_n(z)\}_{n=0}^{\infty}$ be  the sequence of orthonormal polynomials associated with the inner product induced by $\mathbf{M}$.  If  $\mathcal{D}$ is  bounded in $\mathbb{P}[z]$  then the set of zeros of the orthonormal polynomials associated to $\mathbf{M}$ is bounded. Moreover,
$$
Z=
 \{ z \in \CC: \text{ there exists } n\in \NN, \varphi_n(z)=0\}\subset \{ z\in \CC  : \;  \vert z \vert < \Vert \mathcal{D} \Vert \}.
$$
\end{prop}

\begin{proof} Assume that  $z_0\in Z$, i.e.   there exists  $n\in \NN$ such that  $\varphi_n(z)=0$. Then, $\varphi_n(z)=(z-z_0)P_{n-1}(z)$ with  $0\neq P_{n-1}\in \mathbb{P}_{n-1}[z]$. Therefore,  $\varphi_n(z)=zP_{n-1}(z)-z_0 P_{n-1}(z)$ and  since  $P_{n-1}$ and  $\varphi_n(z)$ are orthogonal with respect to $\mathbf{M}$ it follows that 

$$\Vert zP_{n-1}(z) \Vert^2 = < \varphi_n(z) + z_0 P_{n-1}(z), \varphi_n(z) + z_0 P_{n-1}(z)>= \Vert \varphi_n\Vert^2 + \vert z_0 \vert^{2}\Vert P_{n-1}\Vert^2$$ 

\noindent therefore
$$
\vert z_0 \vert^2 \Vert P_{n-1} \Vert^2 = \Vert zP_{n-1} \Vert^2 - \Vert \varphi_n \Vert^2 < \Vert z P_{n-1}(z) \Vert^2 = \Vert \mathcal{D}(P_{n-1}) \Vert^2 \leq \Vert \mathcal{D} \Vert^2 \Vert P_{n-1} \Vert^2
$$

\noindent Since  $\Vert P_{n-1} \Vert \neq 0$ it follows that 
$$
\vert z_0  \vert^2 <  \Vert \mathcal{D} \Vert^2
$$

\noindent as we required.

\end{proof}

\noindent As a consequence of the above result, in order to prove that zeros of orthogonal polynomials associated to measures, Sobolev inner products or more generally  matrix Sobolev inner products,  are uniformly bounded we require  conditions in terms of $\mathbf{M}=(c_{i,j})_{i,j=0}^{\infty}$ to assure that the operator $\mathcal{D}$ is bounded.  For standard orthogonal polynomials associated to a measure in \cite{EGT6}  the following  characterization of boundedness of multiplication operator is provided in terms  of the moments: $\mathcal{D}$ is  bounded if and only if  $\displaystyle{\lim_{n\to \infty} \dfrac{c_{n+1,n+1}}{c_{n,n}} <\infty}$. In case of Sobolev orthogonal polynomials, or more generally orthogonal polynomials associated to an HPD matrix, this limit may not exist; nevertheless a certain   necessary condition remains true:

\begin{prop} Let $\mathbf{M}=(c_{i,j})_{i,j=0}^{\infty}$ be an HDP matrix and let $\mathcal{D}$ be the multiplication operator on $\mathbb{P}[z]$.  Assume that   $\mathcal{D}$ is bounded in $\mathbb{P}[z]$ with respect to $<\cdot,\cdot>_{\mathbf{M}}$,
then
$$
\limsup_{n} \dfrac{c_{n+1,n+1}}{c_{n,n}} \leq \Vert \mathcal{D} \Vert^{2}
$$

\end{prop}
\begin{proof} Assume that $\mathcal{D}$ is bounded in $\mathbb{P}[z]$ and consider the monomials $\dfrac{z^n}{\Vert z^n \Vert}$ for all $n\geq 1$. Then, 
$$
\Vert \mathcal{D} \Vert^2 \geq \Vert \mathcal{D}\left(\dfrac{z^n}{\Vert z^n\Vert}\right)\Vert^2=
\Vert \dfrac{z^{n+1}}{\Vert z^n \Vert}\Vert^2 = \dfrac{<z^{n+1},z^{n+1}>}{<z^n,z^n>}=\dfrac{c_{n+1,n+1}}{c_{n,n}}
$$

\noindent Consequently, 
$$
\limsup_{n\to \infty} \dfrac{c_{n+1,n+1}}{c_{n,n}} \leq \Vert \mathcal{D} \Vert^2
$$

\end{proof} 

\begin{corollary} Let  $\mathbf{M}=(c_{i,j})_{i,j=0}^{\infty}$ be an HDP matrix and let $\mathcal{D}$ be the multiplication operator.  Assume that 
$$
\limsup_{n} \dfrac{c_{n+1,n+1}}{c_{n,n}} =\infty
$$
\noindent Then, $\mathcal{D}$ is not bounded  on $\mathbb{P}[z]$. 

\end{corollary}
\begin{remark}  In the general framework of infinite  HPD, even for Sobolev moment matrices,  we do not know if boundedness of the sequence  $\{\frac{c_{n+1,n+1}}{c_{n,n}}\}_{n=0}^{\infty}$ implies boundedness of  associated multiplication operator as in the standard case. 
\end{remark}

\section{Sequentially dominated HPSD matrices and   boundedness of the multiplication operator}

\noindent In \cite{LP} they say that a set of $k+1$ measures $\{\mu_j\}_{j=0}^{k}$ in the complex plane is sequentially dominated 
if   for every $j=1,\dots,k$ ${\it supp}(\mu_j) \subset {\it supp}(\mu_{j-1})$ and
$$
d\mu_j=f_{j-1}d\mu_{j-1}, \qquad \qquad f_{j-1} \in L_{\infty}(\mu_{j-1}),\qquad j=1,\dots, k.
$$

\noindent Later on, in \cite{APRR} they generalize this notion and they say that a vectorial measure $\mu=(\mu_0,\mu_1,\dots,\mu_k)$ in the complex plane of finite positive Borel measures is extended sequentially dominated, (i.e. $\mu \in ESD$) if there exists a constant $C>0$ such that for every $j=1,\dots, k$, 
$$
\mu_j \leq C \mu_{j-1}.
$$
\noindent   A vectorial measure $\mu$ is sequentially dominated if and only if $\mu \in ESD$ and $\mu_0$ is infinitely supported. 

\noindent We here introduce a definition which is a generalization of the above notion for measures with finite moments. In order to do it we need the following definition for general HPSD matrices:

\begin{definition} Let $\{ \mathbf{M}_j\}_{j=0}^{k}$ be
a  set of $k+1$ HPSD matrices. We say that  $\{ \mathbf{M}_j\}_{j=0}^{k}$ is a set of {\it sequentially dominated} SHPD matrices if there exists a positive constant $C>0$ such that 
for every $j=1,\dots, k$ such that
$$
\mathbf{M}_{j} \leq C\mathbf{M}_{j-1}
$$

\end{definition}

\begin{definition} A vectorial measure $\mu=(\mu_0,\mu_1,\dots,\mu_k)$ in the complex plane of  positive Borel measures with finite moments is
{\it matrix sequentially dominated} (or $\mu\in MSD$) if the corresponding sequence of associated moment matrices is sequentially dominated, i.e., there exists a constant $C>0$ such that for every $j=1,\dots,k$
$$
\mathbf{M}(\mu_j) \leq C \mathbf{M}(\mu_{j-1})
$$
\end{definition}

\begin{lem} Let $\mu=(\mu_0,\mu_1,\dots,\mu_k)\in ESD$  and such that $\mu_j$ has finite moments for $j=1,\dots,k$. Then $\mu$ is matrix sequentially dominated, i.e, $\mu \in MSD$. 
\end{lem}

\begin{proof} Assume that there exists a constant such $
\mu_j \leq C \mu_{j-1}  
$ for all $j=1,\dots,k$. 
Since each measure $\mu_j$ has finite moments then  every polynomial is integrable and since 
$
\mu_j \leq C \mu_{j-1}. 
$ then for every polynomial $p(z)\in \mathbb{P}[z]$ we have 
$$
\int \vert p(z) \vert^{2} d\mu_{j} \leq C \int \vert p(z) \vert^{2}d\mu_{j-1}
$$

\noindent Therefore, via the identification in (\ref{aste2}) it holds that 
for every $v\in c_{00}$ 
$$
v\mathbf{M}(\mu_j)v^{*} \leq C v \mathbf{M}(\mu_{j-1})v^{*}. 
$$
\noindent Thus,  for every $j=1,\dots, k$ 
$$
\mathbf{M}(\mu_{j}) \leq C \mathbf{M}(\mu_{j-1}). 
$$
\end{proof}

\noindent  We point out that the converse of Lemma 1  is not true, that is, if a vectorial measure belongs to  MSD then 
 it has not to be extended sequentially dominated. Indeed, consider the following example: 

\begin{ejem} There is a vectorial measure $\mu=(\mu_0,\mu_1)\in MSD$ such that $\mu \notin ESD$. 
Let  $\mu_0={\bf m}$ the Lebesgue measure in the unit circle and let $\mu_1$ be the Lebesgue measure  in the circle $\{z\in \CC: \Vert z \vert =\frac{1}{2}\}$. 
Therefore $\mathbf{M}(\mu_0)=\mathbf{I}$ and $\mathbf{M}(\mu_1)$
is an infinite  diagonal matrix with diagonal entries  $(2^{-2n})_{n=0}^{\infty}$ (see e.g. \cite{EGT6} for details). It is clear that for every $v=(v_0,v_1,\dots,v_n,0,0,\dots)\in c_{00}$, 

$$
v\mathbf{M}_1v^{*}=   \sum_{k=0}^{n} \vert v_k \vert^{2}\dfrac{1}{2^{2n}}\leq 
\sum_{k=0}^{n} \vert v_k \vert^{2} \leq v\mathbf{M}_0v^{*}
$$

\noindent Consequently $\mu \in MDS$ Nevertheless, the vectorial measure $\mu=(\mu_0,\mu_1)$ does not verify 
that there is a constant $C>0$ such that $\mu_1 \leq C\mu_0$. Note that it this was the case then the support of $\mu_1$ would be included in the support of $\mu_0$.  
\end{ejem} 
\bigskip



\noindent We now prove our main result in this section concerning sequentially dominated matrices which is a generalization of the corresponding results for sequentially dominated measures given in \cite{LP}, \cite{LPP}. In order to do it we introduce certain inner products which are more general than Sobolev inner products. For our results is irrelevant the existence of measures.

\begin{definition} Given  a set of  HPSD matrices $\{\mathbf {M}_j\}_{k=0}^{k}$,  with $\mathbf{M}_0>0$ we define the matrix Sobolev inner product on $\mathbb{P}[z]$ as

$$<p(z),q(z)>_{\mathbf{M}_{S}}=\sum_{j=0}^{k}  <p^{(j)},q^{(j)}>_{\mathbf{M}_j}=v\sum_{j=0}^{k} \mathbf{\Lambda}^{j}\mathbf{M}_j\left(\mathbf{\Lambda}^{j}\right)^{*} w^{*}$$   
where ${\mathbf{M}_{S}}=\sum_{j=0}^{k} \mathbf{\Lambda}^{j}\mathbf{M}_j\left(\mathbf{\Lambda}^{j}\right)^{*} $ and $\mathbf{{\Lambda}}^{j}$ as in (\ref{Filk}).


\end{definition}

\begin{remark}
 We consider $\mathbf{M}_0>0$ in the above definition in order to assure that the expression is an inner product well defined in $\mathbb{P}[z]$.
\end{remark} 

\begin{theorem} Let $\{ \mathbf{M}_j\}_{j=0}^{k}$ be a set of sequentially dominated
HPD matrices.  Assume that for every $j=0,\dots, k$ the  associated multiplication operator $\mathcal{D}_{j}$ is bounded in $(\mathbb{P}[z],<\cdot, \cdot>_{\mathbf{M}_j})$. Consider the matrix Sobolev inner product:
$$
<p(z),q(z)>_{\mathbf{M}_{S}}= \sum_{j=0}^{k} <p^{(j)},q^{(j)}>_{\mathbf{M}_j}
$$

\noindent Then the multiplication operator is bounded in $(\mathbb{P}[z],<\cdot,\cdot>_{\mathbf{M}_{S}})$. 

\end{theorem}

\begin{proof} Let $p(z)\in \mathbb{P}[z]$, since $\left(zp(z)\right)^{j-1)}=jp^{(j)}(z)+zp^{j)}(z)$ for every $j=1,\dots,k$ we have that:
$$
<\left(zp(z)\right)^{j)},\left(zp(z)\right)^{j)}>_{\mathbf{M}_j}=<jp^{j-1)}(z)+zp^{j)}(z),jp^{j-1)}(z)+zp^{j)}(z)>_{\mathbf{M}_
j}
$$
$$
=j^2<p^{j-1)}(z),p^{j-1)}(z)>_{\mathbf{M}_j}+<zp^{j)}(z),zp^{j)}(z)>_{\mathbf{M}_j}+2j\mathfrak{R}\left(<
p^{j-1)}(z),zp^{j)}(z)>_{\mathbf{M}_j}
\right) \leq
$$

$$
j^2C_{j}\Vert p^{j-1)}(z) \Vert^2_{\mathbf{M}_{j-1}}+\Vert \mathcal{D}_j \Vert^2 \Vert p^{j)}(z)\Vert^2_{\mathbf{M}_j}+2j\mathfrak{R}\left(<
p^{j-1)}(z),zp^{j)}(z)>\right)_{\mathbf{M}_j} \leq
$$

$$
(j^2C_j+\Vert \mathcal{D}_j \Vert^2 )\Vert p(z) \Vert^{2}_{\mathbf{S}}+2j\mathfrak{R} left(<p^{j-1)}(z),zp^{j)}(z)>_{\mathbf{M}_j}
$$

\noindent By using the Cauchy-Schwarz property of the inner product $<\cdot,\cdot>_{\mathbf{M}_j}$ we have

$$
<p(z),q(z)>_{\mathbf{M}_j}\leq \Vert p(z) \Vert_{\mathbf{M}_j} \Vert q(z) \Vert
_{\mathbf{M}_j}
$$

\noindent Therefore,

$$
\mathfrak{R}\left(<
p^{j-1)}(z),zp^{j)}(z)>_{\mathbf{M}_j}\right)\leq
\vert <p^{j-1)}(z),zp^{j)}(z)>_{\mathbf{M}_j} \vert
\leq \Vert p^{j-1)}(z)\Vert_{\mathbf{M_j}} \Vert \mathcal{D}_{j}(p^{j}(z))\Vert_{\mathbf{M_j}}
$$
$$
\leq
 C_j^{1/2} \Vert \mathcal{D}_j \Vert^2 \Vert p(z) \Vert^2_{\mathbf{S}}
$$

\noindent Therefore for all $j=1,\dots,k$,
$$
<\left(zp(z)\right)^{j)},\left(zp(z)\right)^{j)}>_{\mathbf{M}_j} \leq (j^2C_{j}+\Vert \mathcal{D}_j \Vert +2jC_{j}^{1/2}\Vert \mathcal{D}_j \Vert
)\Vert p(z) \Vert^{2}_{\mathbf{S}}
$$

\noindent And consequently for every $p(z)\in \mathbb{P}[z]$ it holds:
$$
\Vert zp(z)\Vert^{2}_{\mathbf{S}} \leq \left(\Vert \mathcal{D}_0 \Vert+\sum_{j=1}^{k} j^2C_{j}+\Vert \mathcal{D}_j \Vert^2 +2jC_{j}^{1/2}\alpha_j)\right)\Vert p(z)\Vert^2_{\mathbf{S}}
$$

\noindent what means that the operator multiplication is bounded in the space on $\mathbb{P}[z]$ with the matrix Sobolev inner product as we required. 
\end{proof}

\begin{corollary} Let $\{\mathbf{M}_j\}_{j=0}^{k}$ be a set of sequentially dominated
HDP matrices. The  
 set of zeros of the orthogonal Sobolev polynomials associated with the matrix Sobolev inner product induced by $\{\mathbf{M}_j\}_{j=0}^{k}$  is uniformly bounded.
\end{corollary}
\bigskip

\noindent In the sequel we introduce comparable matrix Sobolev inner products: 

\begin{definition}  Given two sets of HPD matrices $\{\mathbf{M}_j\}_{j=0}^{k}$ and $\{\mathbf{M}'_{j}\}_{j=0}^{k}$, we say that the associated matrix Sobolev norms are {\it comparable component by component} if $\mathbf{M}_j \asymp \mathbf{M}'_j$
for every $j=0,\dots, k $. 
\end{definition}
\begin{prop} Let $\{\mathbf{M}_j\}_{j=0}^{k}$ a set of HPSD matrices. The following are equivalent: 
\begin{enumerate}
\item $\{\mathbf{M}_j\}_{j=0}^{k}$ is sequentially dominated. 
\item The matrix Sobolev norms   associated to the sets 
 $\{\mathbf{M}_j\}_{j=0}^{k}$ and 
 $\{\mathbf{M}'_j\}_{j=0}^{k}$ where $\mathbf{M}'_{j}=
 \mathbf{M}_{j}+ \mathbf{M}_{j+1}+ \dots +\mathbf{M}_{k}$ are comparable component by component. 
\end{enumerate}
\end{prop}

\begin{proof} Assume that the set of matrices $\{\mathbf{M}_j\}_{j=0}^{k}$ is sequentially dominated then there exists $c>0$ such that for all $j=1,\dots,k$ it holds
$$
\mathbf{M}_j\leq C\mathbf{M}_{j-1}.
$$
\noindent Then, there exists $C>0$ such that for all $j=1,\dots,k$ and $i=0,1,\dots,k-j$
$$
\mathbf{M}_{j+i}\leq C \mathbf{M}_{j}
$$

\noindent Then, for every $j=1,\dots,k$ it holds 
$$
\mathbf{M}_j \leq \mathbf{M}_{j}+M_{j+1} +\dots+ \mathbf{M}_k \leq C(\mathbf{M}_j+\dots+\mathbf{M}_j)=C(k-j+1)\mathbf{M}_j.
$$

\noindent On the other hand, assume that for all $j=1,\dots, k$
$$
 \mathbf{M}_{j}+ \mathbf{M}_{j+1}+ \dots +\mathbf{M}_{k} \asymp \mathbf{M}_{j}
 $$
\noindent Then, there exists $C>0$ such that for all $j=1, \dots,k$
$$
\mathbf{M}_j \leq  \mathbf{M}_{j}+  \dots +\mathbf{M}_{k} \leq C\mathbf{M}_j
$$

\noindent and therefore the set 
 $\{\mathbf{M}_j\}_{j=0}^{k}$ is sequentially dominated. 
\end{proof}

\noindent In  \cite{JR1} it is proved that, for Sobolev inner products  supported in the real line,  the boundedness of the multiplication operator implies that the corresponding  Sobolev norm is essentially sequentially dominanted. Essential sequential domination means that the given Sobolev norm is equivalent to another Sobolev norm which is sequentially dominated. Following basically the same arguments, in \cite{APRR} the authors prove a similar result for measures supported in the complex plane. As a consequence of our results matrix sequentially dominated measures are, in particular, essentially sequentially dominated measures as we prove: 

\begin{corollary}  Let  $\mu=(\mu_0,\mu_1,\dots,\mu_k)$ be a vectorial measure of  infinitely supported measure.  If  $\mu$ is 
{\it matrix sequentially dominated} (or $\mu\in MSD$) then $\mu$ is essentially sequentially dominated. 
\end{corollary}
\begin{proof} If $\mu$ is matrix sequentially dominated then $\{\mathbf{M}(\mu)\}_{j=0}^{k}$ is a set of sequentially dominated matrices. Therefore by 
 Proposition 3 the  matrix Sobolev norms   associated to the sets 
 $\{\mathbf{M}(\mu_j)\}_{j=0}^{k}$ and 
 $\{\mathbf{M}(\mu_j')\}_{j=0}^{k}$ where $\mu'_{j}=
 \mu_{j}+ \mu_{j+1}+ \dots +\mu_{k}$ are comparable component by component. Then, obviously, the associated matrix Sobolev norms are equivalent. This means that $\mu$ is essentially sequentially dominated.  

\end{proof}

\section{Location of zeros of orthogonal polynomials via generalized eigenvalues. } 

\noindent We use the notion of generalized eigenvalues in order to give a characterization of sequentially dominated HPD matrices. We need the following lemma: 

\begin{lem} let $\mathbf{A},\mathbf{B}$ be HPD.   The following are equivalent: 
\begin{enumerate}
    \item There exists $c>0$ such that $\mathbf{A} \leq c\mathbf{B}$
    \item $\beta(\mathbf{A},\mathbf{B})<\infty$
\end{enumerate}
\end{lem}

\begin{proof} Since $\mathbf{B}>0$ for every $v\in c_{00}\setminus \{0\}$ 
$$
v\mathbf{A}v^{*} \leq c  v\mathbf{B}v^{*} \Longleftrightarrow \dfrac{v\mathbf{A}v^{*}}{v\mathbf{B}v^{*}} \leq c
$$
\noindent The results now follows easily from the definition of $\beta(\mathbf{A},\mathbf{B})$.
\end{proof}

\noindent Using the notion of generalized eigenvalues and the above lemma we may give the following characterization of sequentially dominated matrices. 

\begin{prop}  Let $\{ \mathbf{M}_j\}_{j=0}^{k}$ be
HPD matrices  for $j=0,\dots,k$. The following are equivalent:

\begin{enumerate}
\item $\{\mathbf{M}_j\}_{j=0}^{k}$ is a set of {\it sequentially dominated}  matrices. 

\smallskip 

\item $\beta(\mathbf{M}_{j},\mathbf{M}_{j-1})<\infty$ for all $j=1,\dots,k$.
\end{enumerate}
\end{prop}

\noindent Above proposition lets to provide a criteria  of boundedness of the set of zeros of Sobolev polynomials associated to HPD matrices in terms of generalized eigenvalues:  

\begin{corollary} Let $\{ \mathbf{M}_j\}_{j=0}^{k}$ be
HPD matrices  for $j=0,\dots,k$.  
Assume that 
$$\beta(\mathbf{M}_j,\mathbf{M}_{j-1})<\infty \qquad  j=1,\dots,k$$

\noindent Then, the   
 set of zeros of the orthonormal  polynomials associated with the matrix Sobolev inner product $<\cdot,\cdot>_{\mathbf{M}_{S}}$ is  bounded.
\end{corollary}

\noindent An easy consequence of Lemma 1 and Corollary 4 is the following

\begin{corollary} Let $\mu=(\mu_0,\dots,\mu_k)$ be a vectorial measure in MSD with all measures infinitely supported. Then the set of zeros of Sobolev orthogonal polynomials is bounded. 

\end{corollary}

\bigskip

\noindent Using the results in  \cite{EGT6} relating with the generalized eigenvales and  the notion of the polynomially convex hull of the support of a measure $\mu$ (see e.g. \cite{Conway}), denoted by $P_C({\it supp}(\mu))$,  we obtain: 

\bigskip

\begin{corollary} Let $\mu=(\mu_0,\dots,\mu_k)$ be a vectorial measure in MSD with all measures infinitely supported. Then for every $j=1,\dots,k$
$$
{\it supp}(\mu_j) \subset P_{C}({\it supp}(\mu_{j-1})).
$$
\end{corollary}

\begin{remark} It is interesting to point out that for an  ESD vectorial measure $\mu=(\mu_0,\dots,\mu_k)$ it is obvious  that ${\it supp}(\mu_j) \subset {\it supp}(\mu_{j-1})$. This is a more restricted condition for the supports  than  in Corollay 6 for MSD vectorial measures. 
\end{remark}


\noindent In the sequel we analyze the case of sets of two sequentially dominated matrices  when one of them is the identity matrix, that is, the moment matrix associated to the Lebesgue measure  in the unit circle ${\bf m}$. Note that in this case, generalized eigenvalues coincide with standard eigenvalues.

\begin{corollary} Let $\mathbf{M}$ be an HPD matrix  and let $\lambda_n,\beta_n$ the smallest and largest eigenvalue of the finite  sections. Let  $<\cdot,\cdot>_{\mathbf{M}_{S}}$ the matrix inner Sobolev product associated to $\{\mathbf{M}_0,\mathbf{M}_1\}$. In the following two situations 

\begin{enumerate} 
\item  If $\mathbf{M}_0=\mathbf{I}$,  $\mathbf{M}_1=\mathbf{M}$ and $\displaystyle{\lim_{n\to \infty} \beta_n<\infty}$,

\item  If  $\mathbf{M}_0=\mathbf{M}$,  $\mathbf{M}_1 =\mathbf{I}$ and $\displaystyle{\lim_{n\to \infty} \lambda_n>0}$,  
\end{enumerate}

\noindent the set of zeros of Sobolev orthogonal polynomials is bounded. 
\end{corollary} 

\smallskip

\noindent In the case of measures the above corollary can be rewritten as follows: 

\begin{corollary} Let $\mu$ be an infinitely and compactly supported measure  in $\CC$ and let $\lambda_n,\beta_n$ the smallest and largest eigenvalue of the finite sections of the moment matrix   $\mathbf{M}(\mu)$. In the following two cases: 

\begin{enumerate} 
\item  $\displaystyle{\lim_{n\to \infty} \beta_n<\infty}$ and  
$$
< p(z), p(z)>_{\mathbf{S}}= \int \vert p(z) \vert^2 \; d{\bf m}+ \int \vert p'(z)\vert^2 d\mu, 
$$

\item  $\displaystyle{\lim_{n\to \infty} \lambda_n>0}$  and 
$$
< p(z), p(z)>_{\mathbf{S}}=   \int \vert p(z) \vert^2 \; d\mu+ \int \vert p'(z)\vert^2 d{\bf m} 
$$

\end{enumerate}
\noindent the set of zeros of Sobolev orthogonal polynomials is bounded. 

\end{corollary}

\begin{remark}\noindent These conditions in terms of  eigenvalues are sufficient conditions to assure boundedness of zeros of polynomials, thus these results do not provide information in  the case that one of the measures is the Lebesgue measure in the unit circle and the other measure verifies that $\displaystyle{\lim_{n\to \infty} \lambda_n=0}$ and $\displaystyle{\lim_{n\to \infty} \beta_n =\infty}$. We provide an example of this situation:
\end{remark}

\begin{ejem}
Let 
 $\mu$ be the normalized Lebesgue measure in $\{z\in \CC: \vert z-1 \vert =1\}$ and the inner Sobolev product: 
$$
<p(z),q(z)>_{\mathbf{S}}=
\int p(z)\overline{q(z)}d\mu+\int p'(z)\overline{q'(z)}d{\bf m}.
$$

\noindent  The associated Sobolev moment matrix is 
$\mathbf{M}_{\mathbf{S}}=\mathbf{M}+ \mathbf{\Lambda}\mathbf{I}\mathbf{\Lambda}^{*}
$ where $\mathbf{M}=\mathbf{M}(\mu)$, and it can be easily checked that $\mathbf{M}_{\mathbf{S}}=\left ( \binom{i+j}{i} + \delta_{ij} i^2 \right )_{i,j=0}^{\infty}$, i.e.,
 \[ \mathbf{M}_{S}= \left ( \begin {array}{cccccc} 1&1&1&1&1& \ldots\\ \noalign{\medskip}1&3&3&4&5& \ldots
\\ \noalign{\medskip}1&3&10&10&15& \ldots\\ \noalign{\medskip}1&4&10&29&35& \ldots
\\ \noalign{\medskip}1&5&15&35&86& \ldots \\
\vdots & \vdots & \vdots & \vdots  & \ddots
\end {array} \right )=\left ( c_{j,k} \right )_{j,k=0}^{\infty}, \]

\noindent  since the moment matrix $\mathbf{M}_0=\left (\binom{i+j}{i} \right )_{i,j=0}^{\infty}$ (see e.g. \cite{EGT1}) and $\mathbf{M}_1=\mathbf{I}$. In  this case the Cauchy-Schwartz condition $c^2_{n,n} \leq c_{n-1,n-1}c_{n+1,n+1}$ does not hold. Nevertheless, 

 \[ \lim_{n \to \infty} \frac{c_{n+1,n+1}}{c_{n,n}} =
\lim_{n \to \infty}  \dfrac{\dfrac{2^{n+1}}{(n+1)!}+ (n+1)^2}{\dfrac{2^{n}}{(n)!}+ n^2}=1.\]

\noindent Moreover  $\displaystyle{\lim_{n\to \infty} \lambda_n=0}$ (see \cite{EGT1}); therefore  corollary $8$ does not provide any information about the location of zeros of Sobolev polynomials.  
\end{ejem}

\noindent In the particular case of measures  supported in the closed unit disk, using the results in \cite{EGT1} we have 

\begin{corollary} Let $\mu$ be a measure with support in $\overline{\mathbb{D}}$ and such that $\mu/\mathbb{T}= w(\theta)d\theta$. In the following two cases, 
\begin{enumerate} 

\item  ess sup $w(\theta)<\infty$ and 
$$
< p(z), p(z)>_{\mathbf{S}}=  \int \vert p(z) \vert^2 \; d{\bf m}+ \int \vert p'(z)\vert^2 d\mu, \text {or} 
$$

\item ess inf $w(\theta)>0$ and  
$$
< p(z), p(z)>_{\mathbf{S}}=  \int \vert p(z) \vert^2 \; d\mu+ \int \vert p'(z)\vert^2 d{\bf m}, 
$$

\end{enumerate} 
\noindent  the set of zeros of Sobolev orthogonal polynomials is bounded. 
\end{corollary}

\noindent

\begin{remark} This result generalizes the corresponding result of \cite{LPP} for the particular  case of a the vectorial measure $({\bf m},\mu)$ which are sequentially dominated, that is, when the  measure $\mu$ is  supported in $\mathbb{T}$ and verifies $\mu=f(z)d{\bf m}$ for some $f(z)\in L^{\infty}({\bf m})$. By  Szego's results  such measures are measures supported in $\mathbb{T}$ with $\displaystyle{\lim_{n\to \infty} \beta_n<\infty}$. Thus, as a consequence of our results the restriction of being supported in $\mathbb{T}$ is not required; anyway by \cite{EGT1} measures verifying $\displaystyle{\lim_{n\to \infty} \beta_n<\infty}$ must be supported in $\overline{\mathbb{D}}$. 
\end{remark}



\begin{thebibliography}{Man96}


\bibitem{APRR} Alvarez,V.; Pestana,D; Rodr\'{\i}guez,J.M.; Romera,E. \emph{Weighted Sobolev spaces on curves} J. Approx, Theory {\bf 119} (2002) 41-85

\bibitem{Bhatia} Bhatia, R. \emph{Matrix Analysis} Graduate Texts in Matematics; 169 , Springer (1997)


\bibitem{Berg}  Berg, C. and  Szwarc,R.  The Smallest Eigenvalue of Hankel Matrices, Constr. Approx. 34 (2011) 107--133. https://doi.org/10.1007/s00365-010-9109-4.

\bibitem{Berg-Duran}Berg, C and  Dur\'{a}n,A.J.  Analytic function associated to positive definite infinite matrices, J. Math. Anal. Appl. 315 (2006) 54--67.

\bibitem{APJ} D\'{\i}az Gonz\'{a}lez,A.; Pijeira-Cabrera, H. and Quintero-Roba,J. \emph{Polynomials of Least Deviation from Zero in Sobolev $p$-norm} Bull. Malays. Math. Sci. Soc.  (2022) 
\bibitem{Conway} Conway, J.B. \emph{A Curse in Functional Analysis}
       Graduate Texts in Mathematics, 96, Springer-Verlag, New York, 1985.
\bibitem{Duran-Saff} Dur\'{a}n, A. and Salff, E. \emph{Zero location for non standard orthogonal polynomials } J. Approx. Theory {\bf 113} (2001) 127-141 



\bibitem{EGT1} Escribano,C; Gonzalo, R. and Torrano, E. \emph{ Small Eigenvalues of Large Hermitian
moment matrices.} J. Math. Anal. Appl. {\bf 374}  pp. 470-480 (2011)

\bibitem{EGT5} Escribano,C; Gonzalo, R. and Torrano, E. \emph{ A characterization of polynomial density on curves via matrix algebra} Mathematics 7, 1231,  (2019) doi:10.3390/math 712131.
 
\bibitem{EGT6}  Escribano,C; Gonzalo, R. and Torrano, E.  \emph{Smallest and largest generalized eigenvalues of large moment matrices and some applications} J. Math. Anal. App. {\bf 521(2)} (2023)

\bibitem{LP} L\'{o}pez Lagomasino,G. , Pijeira Cabrera, H. \emph{Zero location and n-th root asymptotics of Sobolev orthogonal polynomials} J. Approx. Theory {\bf 99} (1999) pp. 30-43


\bibitem{LPP} L\'{o}pez Lagomasino,G. , Pijeira Cabrera, H. and P\'{e}rez Izquiero,I. \emph{Sobolev orthogonal polynomials in the complex plane} J. Comp. Appl. Math. {\bf 127} (2001) pp. 
 219-230

\bibitem{Finkelshtein} Mart\'{\i}nez Finkelshtein, A. \emph{Analytic aspects of Sobolev orthogonal polynomials revisited} {\bf 127}  Journal of Computational and Applied Mathematics (2001) pp.255-266


\bibitem{JR1} Rodr\'{\i}guez, J.M. \emph{Multiplication operator in Sobolev spaces with respect to measures.} J. Approx. Theory {\bf 109}  (2001) 157-197

\bibitem{JR3} Rodr\'{\i}guez,J.M. \emph{ A simple characterization of weighted Sobolev spaces with bounded multiplication operator} J. Approx. Theory {\bf 153} (2008) 53-72

\end{thebibliography}
\end{document}